\documentclass[sn-nature]{sn-jnl}

\usepackage{amssymb}
\usepackage[T1]{fontenc}
\usepackage{float}
\usepackage{amsmath}
\usepackage{xfrac}
\usepackage{amsfonts}
\usepackage{mathrsfs}
\usepackage{graphicx}
\usepackage{manyfoot}

\newcommand{\mytitle}{Model-Based Reinforcement Learning Control of Reaction-Diffusion Problems}
\newcommand{\myabstract}{Mathematical and computational tools have proven to be reliable in decision-making processes. In recent times, in particular, machine learning-based methods are becoming increasingly popular as advanced support tools. When dealing with control problems, reinforcement learning has been applied to decision-making in several applications, most notably in games. The success of these methods in finding solutions to complex problems motivates the exploration of new areas where they can be employed to overcome current difficulties. In this paper, we explore the use of automatic control strategies to initial boundary value problems in thermal and disease transport. Specifically, in this work, we adapt an existing reinforcement learning algorithm using a stochastic policy gradient method and we introduce two novel reward functions to drive the flow of the transported field. The new model-based framework exploits the interactions between a reaction-diffusion model and the modified agent. The results show that certain controls can be implemented successfully in these applications, although model simplifications had to be assumed.}

\begin{document}

\title{\mytitle}

\author*[1]{Christina Schenk} \email{christina.schenk@imdea.org}
\author[1]{Aditya Vasudevan} \email{adityavv.iitkgp@gmail.com}
\author[1]{Maciej Haranczyk} \email{maciej.haranczyk@imdea.org}
\author[1,2]{Ignacio Romero}\email{ignacio.romero@upm.es}

\affil[1]{Universidad Polit\'ecnica de Madrid,
Jos\'e Guti\'errez Abascal, 2, Madrid 29006,
Spain}

\affil[2]{IMDEA Materials Institute,
	Eric Kandel 2, Tecnogetafe, Madrid 28906,
	Spain}

\abstract{\myabstract}
\maketitle

\section{Introduction}
\label{sec:intro}
Tools of optimal control have found their application in many fields, ranging from financial markets to the food industry and self-driving vehicles \cite{sivamayil_systematic_2023, schenk17, hull2013optimal}. One particularly difficult area is the control of solutions governed by initial boundary value problems. For instance, in mechanical engineering, it can sometimes be convenient to optimize the design of a thermally conductive part to maximize its thermal dissipation while employing as little material as possible. In this example, the model that needs to be controlled is described by the heat equation, the canonical parabolic partial differential equation.

In a seemingly unrelated field, the spread of the SARS-CoV-2 virus --- and similar pandemics --- can be studied with transport models that predict the evolution of the population that is infected by the disease, also described by parabolic partial differential equations. The goal of an optimal control here could be to minimize the number of infected people while, simultaneously, minimizing mobility restrictions (see \cite{napolitano_impact_2022}, which reviewed 17,000 works on computational modeling of pandemic spread). Roughly speaking, the difficulty here lies in controlling the contagion while maximizing the free movement of people.

Most studies of pandemic spread are based on stochastic or deterministic compartment models based on ordinary differential equations. Hence, to study the potential control of the spread, researchers have employed transmission dynamics \cite{shen_mathematical_2021, k_asamoah_mathematical_2020}. Most of these works do not account for the spatial distribution of infections. Instead, they focus on the evolution of lumped quantities (or compartments) such as infected, recovered, and deceased populations. As advanced, modeling the spatial \emph{and} time changes of contagious diseases requires the use of partial differential equations, and their detailed control involves the consideration of the optimal spatial distribution of mobility/diffusivity that modulates transmission dynamics. Remarkably, the most precise models of the aforementioned type are reaction-diffusion initial boundary value problems, such as the heat equation, and have been previously considered in the literature \cite{mammeri_reaction-diffusion_2020, tu_reactiondiffusion_2023, Grave2021AssessingTS, GRAVE2022115541, guglielmi_delay_2022, viguerie_diffusionreaction_2020}.

Accepting that a pandemic space-time evolution can be modeled with a reaction-diffusion partial differential equation, this begs the question of whether this description can facilitate the design of effective control methods. It is well known that, given some initial conditions and spatial distribution of the diffusivities, it is possible to obtain an accurate description of the infections at every point of a certain region and every instant of a given period of time by solving the forward problem. The inverse problem, in contrast, is not well-posed even for constant parameters and a vanishing reaction term, hinting at the fact that the control of the problem solution is hard.

One powerful way to solve control problems, in general, is reinforcement learning (RL). In this evolving field of machine learning, an agent is trained by its interaction with an environment to maximize a reward. Reinforcement learning has been applied for various applications, from defense/military scenarios to autonomous vehicles, networking, robotics grid world, finance, marketing, natural language processing, Internet of Things security, energy management, recommendation systems, and -- especially -- games \cite{wells_explainable_2021, sivamayil_systematic_2023}. The ability of RL to control such diverse problems suggests that these techniques might also be employed to control transport problems.

Roughly speaking, the applicability of RL algorithms to control problems is limited due to the curse of combinatorial explosion. In real problems, the size of the state-action space is too large, and there are far too many possible combinations of states and actions. This is true, in particular, for the application of RL methods to the control of initial boundary value problems, where the state space is infinite-dimensional. An obvious approach consists of transforming these problems into \emph{approximate}, finite-dimensional, ones. If the space and time granularity of these approximations is fine enough, we might expect their solution to be close to the exact ones and amenable to automatic control. Since the state- and state-action space grow rapidly for real-world applications, most recent efforts have aimed at optimizing existing algorithms applied to toy examples or small-scale case studies \cite{wells_explainable_2021}. Lately, new avenues of research aim to make RL strategies more explainable through human collaboration, visualization, policy summarization, query-based explanations, and verification \cite{wells_explainable_2021}.

Reinforcement learning (RL) has found application in the realm of optimal control and containment strategies to manage COVID-19, particularly in the refinement of epidemic control measures and lockdown policies \cite{doi:10.1080/08839514.2022.2031821,arango2020,PADMANABHAN2021102676,ohi_exploring_2020,khadilkar_optimising_2020}. These studies consider complex epidemiological models. However, to our knowledge, the consideration of spatial dynamics in these RL approaches has yet to be explored.

In numerous applications, the optimization objective extends beyond a single goal, encompassing multiple objectives. As a result, the identification of the Pareto optimal solution becomes pivotal. However, obtaining the Pareto front directly proves to be a challenging task, prompting numerous scholars to explore various approaches in the realm of multi-objective reinforcement learning \cite{hayes_practical_2022, zhang_multi-objective_2023}. These approaches involve formulating a multi-objective Markov decision process, which estimates the values of each considered policy, departing from the conventional approach of estimating a single value for a single policy. As an alternative to the multi-objective Markov decision process formulation approach, a single-weighted policy can be defined. Each of these approaches carries its own set of advantages and disadvantages \cite{hayes_practical_2022}. Given that multi-objective RL algorithms are not implemented in the selected Tensorforce library \cite{tensorforce}, we opt for the single-weighted policy approach in this paper.

In this work, we present a multi-objective RL algorithm, i.e. a model-based stochastic policy-gradient algorithm that is implemented by a modification of the Tensorforce agent and the model-described environment. In addition, we explore its applicability to the control of parabolic diffusion problems and, more specifically, those arising in thermal and disease transport. In the case of thermal problems, the RL algorithm is applied to a control problem maximizing the thermal dissipation of a domain while keeping the temperature bounded; for infection transmission problems, the goal of the proposed RL algorithms will be to maximize mobility, while minimizing the infected population. In both situations, the environment is modeled by a single field or a reaction-diffusion model. The environment and agent interact with one another: the agent sends the control action to the environment and the environment returns the resulting observed state back. In our case, we have employed two connected software tools, one C++-based code to solve the transport problems and one Python-based code based on Tensorforce that is responsible for the RL. In addition to the states, the $L^2$ norms of diffusivities and states for calculating the reward functions are sent from the environment to the agent.

The remainder of this article has the following structure. In~Section~\ref{sec:meth}, the underlying methods are introduced. First, the model, then the optimal control problem with the considered reward functions, and the RL algorithm are introduced. Then, the software setup is described.  In Section~\ref{sec:results}, the proposed methods are applied to two examples, a thermal problem and a disease transmission dynamics problem. The results of considering two reward functions are compared and discussed. The article concludes with a summary of the main findings in Section~\ref{sec:concl}.

\section{Methods}
\label{sec:meth}
\subsection{Model}
\label{sec:models}
The work described in this work is based on models described by the following initial boundary value problem. Let $\Omega$ be a bounded open set in $\mathbb{R}^2$ with smooth boundary $\Gamma=\Gamma_g\cup\Gamma_h$ and $\Gamma_g\cap\Gamma_h=\emptyset$. To study
the evolution during a period of length $T$ of $c:\Omega\times[0,T]\to[0,1]$, the ratio of the infected population, we will use the single field model
\begin{subequations}
	\begin{align}
	&\dot{c}(x,t) = \nabla \cdot (\rho\,\kappa(x)\,\nabla c(x,t)) + \mathcal{R}(c)\ ,
 &(x,t)\in\Omega\times(0,T) \label{eq:model1}\\
 &c(x,t) = g(x,t)\ ,  &(x,t)\in\Gamma_g\times(0,T)\\
    &\kappa(x)\nabla c(x,t)\cdot\hat{n}(x,t) = h(x,t)\ ,
    &(x,t)\in\Gamma_h\times(0,T)\\
	&c(x,0) = c_0(x) & x\in\Omega   \ .
	\end{align}
\label{eq:model}
\end{subequations}
Here, $\rho$ is the population density, $\kappa$ the diffusivity field, $h$ the influx of infected population through the boundary $\Gamma_h$, $g$ the known infected population ratio on $\Gamma_g$, and $c_0$ the initial value of $c$. The symbols $\nabla, \nabla\cdot$ denote the gradient and divergence operators; the notation $\dot{c}$ indicates the partial derivative with respect to time.

The term $ \mathcal{R}(c)$ is the \emph{reaction} of the model and here it is selected as
\begin{equation}
\mathcal{R}(c) = (\beta - \gamma)c - \beta c^2.
\label{eq:R}
\end{equation}
This choice corresponds to a SIS (Susceptible - Infectious - Susceptible)  model, the most basic compartmental model for disease transport considering only two populations, the susceptibles $S$ and the infectious $I$ \cite{Kuhl}. This two-compartment model can be reduced to a single field, where the infection compartment $I = c$ and the uninfected or the susceptible compartment is $S = 1 - c$. The parameter $\beta$ is the contact rate and $\gamma$ is the infectious rate of the population. If $\sfrac{\beta}{\gamma} > 1$, then this results in an endemic equilibrium with $I_{t\to \infty} = 1 - \sfrac{\gamma}{\beta}$ and $S_{t \to \infty} = \sfrac{\gamma}{\beta}$ \cite{Kuhl}. Using this formulation, the infected population ratio must assume values in $[0,1]$. In this work, we fix $\gamma = 1$ and typically choose $\beta >1$.

Remarkably, the initial boundary value problem~\eqref{eq:model} also describes the thermal evolution on a plane domain $\Omega$ if $c$ is interpreted as the temperature, $\kappa$ the thermal diffusivity, $\rho$ the density, $h$ the known value of the supplied heat flux through $\Gamma_h$, $g$ a known thermal field, and $c_0$, the initial distribution of the temperature. For the thermal problem, the reaction term can also be chosen to be nonlinear as in eq. (\ref{eq:R}), with two reaction constants $\beta,\gamma$. For simplicity, for the remainder of this work, we refer to $c$ as infections although, as indicated, it refers to the temperature in a thermal problem.

The identification of the disease and thermal problem opens the door to using the numerous tools available for thermal analyses and, in particular, finite element codes. On the flip side, model~\eqref{eq:model} might represent disease spread too poorly, given that a very simplistic epidemiological model is considered.

\subsection{Space and time numerical approximation}
\label{sec:finite}
To approximate the solution of the initial boundary value problem \eqref{eq:model}, we study its variational form as the basis of the finite element formulation and a time discretization based on an implicit Euler method. To introduce it, let $H^1_g(\Omega)$ denote, as customarily, the Hilbert space of square-integrable functions, with square-integrable (first) derivative, that vanishes on $\Gamma_g$. Multiplying the balance equation~\eqref{eq:model1} by an arbitrary function $\delta c\in H^1_g(\Omega)$, integrating by parts, and making use of the boundary conditions, it follows that
\begin{equation}
    \begin{aligned}
   \label{eq-weak}
	\int_{\Omega}\Big(\dot{c}(x,t)\;\delta c(x) + &\rho\,\kappa(x)\,\nabla c(x,t) \cdot \nabla \delta c(x)
    + \mathcal{R}(c)\delta c(x) \Big) dV \\
    &= \int_{\Gamma_h} h(x,t)\,\delta c(x)\; dV
    \end{aligned}
\end{equation}
holds for all $\delta c$. To discretize in time this variational equation, we will use the implicit Euler method. For that, let $0=t_0<t_1<\cdots<t_n=T$ be $n+1$ time instants and let $\Delta t_i=t_{i+1}-t_i$. Then, we define $c_i(x)$ to be an approximation to $c(x,t_i)$ and replace eq. (\ref{eq-weak}) with the semi-discrete equation
\begin{equation}
\begin{aligned}
   \label{eq-semidiscrete}
	\int_{\Omega}
 \Big(
 \frac{c_{i+1}(x)-c_i(x)}{\Delta t_i}
 &\;\delta c(x) +
 \rho\,\kappa(x)\,\nabla c_{i+1}(x) \cdot \nabla \delta c(x)\\
    & + \mathcal{R} (c_{i+1}(x))\delta c(x)\Big) dV
    = \int_{\Gamma_h} h_{i+1}(x)\,\delta c(x)\; dA\ ,
\end{aligned}
\end{equation}
with $h_{i+1}(x)=h(x,t_{i+1})$. To complete the discretization of the problem, we use finite elements in space. For that, consider a triangulation of $\Omega$ into $n_{el}$ elements connected at $n_n$ nodes defining a set  $\{N_a(x)\}_{a=1}^{n_n}$ of piecewice linear shape functions with $N_a:\Omega \to[0,1]$. These functions define a finite-dimensional space
\begin{equation}
    V^h = \{ u^h\in H^1_g(\Omega),\ u^h(x) = \sum_{a=1}^{n_n} N^a(x)\,u^a\ ,
    u^a\in\mathbb{R}\} \subset H^1_g(\Omega)\ .
\end{equation}
Then, given an arbitrary function $g^h_{i}(x)=\sum_{a=1}^{n_n} N^a(x)\,g^a_{i}$, not in $H^1_g(\Omega)$ but with value $g^h_i(x)=g(x,t_{i})$ on $\Gamma_g$, the finite element solution and weight functions are defined, respectively, as
\begin{equation}
\label{eq-expansions}
	c^h_{i}(x) = g^h_i(x) + \sum_{a=1}^{n_h} N^a(x)\,c^a_i,
 \qquad
 \delta c^h(x) = \sum_{a=1}^{n_h} N_a(x)\,\delta c_a\ .
\end{equation}
The final, fully discrete, approximation of the initial boundary value problem eq. (\ref{eq:model})
is obtained by replacing $c_i$ by $c^h_i$ and $\delta c$ by $\delta c^h$ in eq. (\ref{eq-semidiscrete}). Given $c_j^h$ for $j=0,1,\ldots,i$ the approximate solution at time $t_{i+1}$ is the function $c_{i+1}^h$ of the form~\eqref{eq-expansions} that satisfies
\begin{equation}
    \begin{aligned}
\int_{\Omega}
 \Big(
 \frac{c^h_{i+1}(x)-c^h_i(x)}{\Delta t_i}
 \;\delta c^h(x) &+
 \rho\,\kappa(x)\,\nabla c^h_{i+1}(x) \cdot \nabla \delta c^h(x)\\
    & + \mathcal{R} (c_{i+1}^h(x))\;\delta c^h(x)\Big) dV
    = \int_{\Gamma_h} h_{i+1}(x)\,\delta c^h(x)\; dA\
    \end{aligned}
\label{eq:weakfe}
\end{equation}
for all $\delta c^h\in V^h$.

Notably, the space and time discrete problem of evolution described by the variational equation~\eqref{eq:weakfe} can be recast
as an optimization problem. To see this, let us assume that the solution $c_i^h$
at time $t_i$ is known and consider the incremental potential
\begin{equation}
\begin{aligned}
    I[c^h_{i+1}] =
\int_{\Omega}
 &\left(
 \frac{|c^h_{i+1}(x)-c^h_i(x)|^2}{2\,\Delta t_i}
 +
 \rho\,\kappa(x) \frac{|\nabla c^h_{i+1}(x)|^2}{2}
   +
   \mathcal{H}(c_{i+1}^h(x)) \right)\; dV\\
  &- \int_{\Gamma_h}
 h_{i+1}(x)\,c^h_{i+1}(x)
 \; dA\ ,
 \end{aligned}
\end{equation}
that employs
\begin{equation}
  \label{eq-hh}
  \mathcal{H}(c_{i+1}^h(x)) = (\beta - \gamma) \frac{|c^h_{i+1}(x)|^2}{2} - \beta \frac{|(c^h_{i+1}(x))|^{3}}{3}\ .
\end{equation}
It is straightforward to verify that the stationarity condition of this functional coincides with eq. (\ref{eq:weakfe}), a property that can be exploited to develop optimal control
strategies.

\subsection{Control}
\label{sec:cont}
\begin{figure}[H]
	\centering
	\includegraphics[scale=.7]{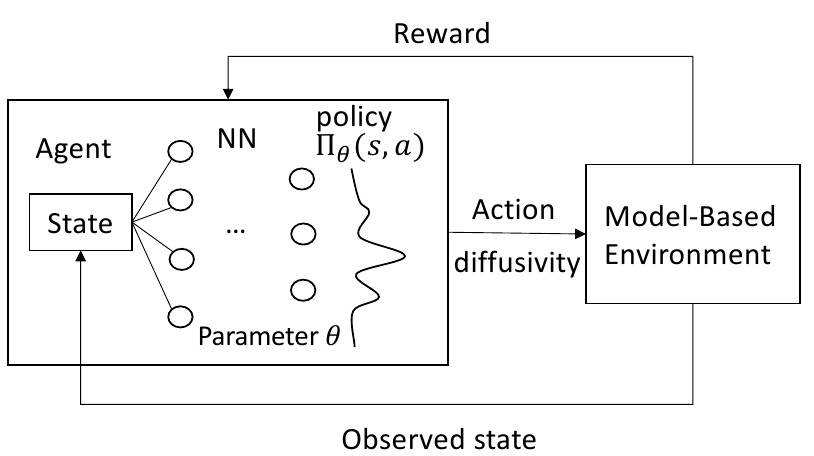}
	\caption{Overview of the RL algorithm adapted from \cite{mao_resource_2016}.}
	\label{Fig:RLAlgOverview}
\end{figure}

To provide decision support and to control advection-reaction-diffusion-type problems (e.g. disease transmission or thermal dynamics), we have developed a model-based reinforcement learning control framework. This framework (cf. fig.~\ref{Fig:RLAlgOverview}) relies on a stochastic policy gradient descent algorithm. For the sake of clarity, we provide a concise overview of how this framework operates within our specific context. Readers seeking a more in-depth understanding of reinforcement learning concepts can refer, for example, to \cite{Sutton1998}.

In this setup, an agent, typically a neural network, interacts with its environment, which here is represented as a model (the finite element simulation of the problem as defined in eq. (\ref{eq:model}). The agent observes the state function $c$ (e.g., infections or temperature) and takes actions on the control variable $\kappa$, modifying its spatial distribution, based on the current observed state $c$. Each action on $\kappa$ results in a reward for the agent. Since the goal of the agent is to maximize this reward, depending on the definition of the latter, the agent is provided with feedback on the effectiveness of its actions. Typically, the RL agent is trained for a maximum number of episodes, where each episode contains a sequence of states, actions, and rewards ending with the terminal state.

Given a reward function and some state $c$, the behavior of an RL agent is completely defined by its \emph{policy}~$\pi$, a stochastic map from the state space to the set of possible actions over the control variable $\kappa$. This mapping is probabilistic since it must balance two strategies: \emph{exploiting} the knowledge of past behaviors and \emph{exploring} new avenues that might lead to larger rewards. The effectiveness of such a delicate balance determines the efficiency of the agent, but more importantly its chances of successfully controlling the environment. The policy is thus dependent on the reward function, and its selection is nontrivial because the mathematical description of a model (e.g., eq. (\ref{eq:model}) is agnostic about its desirable evolution. This choice is incumbent on the RL control and might have a non-negligible effect on the performance and stability of the solution strategy. In the case at hand, the guiding idea is that the agent should strive to keep infection ratios as small as possible while keeping diffusivities large (so that the population is allowed to move freely). The RL agent must balance these mutually opposed goals and attempt to reach a sort of Pareto optimal which is determined by the designer of the policy and reward function.

We will test two reward functions leading to RL algorithms that promote high diffusivities and discourage infections. Vaguely speaking, the first reward is selected to favor the global increase of diffusivities while keeping the infections from increasing; the second one, in turn, is chosen to reduce the infection rate while keeping the diffusivities from decreasing. In more precise terms, the infections and the diffusivities are fields defined on the domain $\Omega$ and the reward functions employ their $L^2$ norm to measure their ``size''.

Using the notation introduced in Sections~\ref{sec:models} and \ref{sec:finite}, the two reward functions mentioned before can be formulated as:
\begin{equation}
	R_{diff}(\kappa_i, c_i, \kappa_0, c_0, c_0^{bef}) =
	\omega_1\frac{\|\kappa^l_i\|}{\|\kappa^0_0\|}-\omega_2 \max\left(0, \frac{\|c_i^l\|-\|c^{bef}_i\|}{\|c^0_0\|}\right),
	\label{eq:R1}
\end{equation}
\begin{equation}
	R_{state}(\kappa_i, c_i, \kappa_0, c_0, \kappa^{bef}_0) =
	    -\omega_3\frac{\|c_i^l\|}{\|c_0^0\|}+\omega_4 \min\left(0, \frac{\|\kappa^l_i\|-\|\kappa_i^{bef}\|}{\|\kappa_0^0\|}\right)\ ,
	\label{eq:R2}
\end{equation}
for $i=0,\ldots,n$ time instants, $l=0,\ldots, n_{ep}-1$ episodes, and $\|\cdot\|$ representing the $L^2$ norm of fields defined on $\Omega$. In these expressions, $c_i^l, \kappa_i^l$ denote, respectively, the diffusivity at time instant $i$ and episode $l$, respectively. The optimization weights for the four terms are defined by $\omega_m, m=1,\ldots 4$.

The goal of an RL algorithm is to maximize its reward. Using an RL method with reward~\eqref{eq:R1}, we aim to maximize the diffusivities while penalizing an infection increase compared to before training. Similarly, an RL aims to minimize infections while penalizing the reduction of mobility when it is based on reward~\eqref{eq:R2}. In detail, the second term of eq. (\ref{eq:R1}) compares the (norms of the) current infections with the ones from before training, denoted by $c_i^{bef}$. The second term of eq. (\ref{eq:R2}) compares the (norms of the) current diffusivities with those from before training, the latter being denoted as $\kappa_i^{bef}$. As the diffusivities and infections have different dimensions, all terms are scaled by the norms of the initial values.

%def reward_compute_norm_cont_weight_normalized_constr(self):
%if self.thetanorm_before == 0.0:
%rew = self.current_kappa_norm[0]/self.current_kappa_norm0[0]
%else:
%rew = self.current_kappa_norm[0]/self.current_kappa_norm0[0] - 100*(max(0.0, (self.current_theta_norm[0]-self.thetanorm_before[self.counter])/self.current_theta_norm0[0]))
%return rew
%
%
%def reward_compute_norm_cont_weight_normalized_constr_theta(self):
%if self.kappanorm_before == 0.0:
%rew = -self.current_theta_norm[0]/self.current_theta_norm0[0]
%else:
%rew = -self.current_theta_norm[0]/self.current_theta_norm0[0] + 100*(min(0.0, (self.current_kappa_norm[0]-self.kappanorm_before[self.counter])/self.current_kappa_norm0[0]))
%return rew
\subsubsection{Algorithm}
\label{sec:alg}
To address the control problem outlined in this study, we developed a custom reinforcement learning algorithm in Python, leveraging the Tensorforce \cite{tensorforce} agent. The implementation details are elucidated in the subsequent paragraphs. Our approach adopts a stochastic policy-gradient-based strategy, utilizing the Adam optimizer, a stochastic gradient descent method that estimates first- and second-order moments adaptively \cite{Kingma2015Adam}. Furthermore, we incorporated a multi-step meta-optimizer that employs the Adam optimizer iteratively for 10 times with 10 line search iterations. Gradient updates occur after reaching the maximum number of episode time steps, with a selected learning rate of $8\times 10^{-5}$.

The exploration parameter is set to $0.1$, and a dense neural network with an activation function $\tanh$ and a layer size of 128 is used. Our chosen objective aligns with the policy gradient type. For the reward estimation, we opt for a horizon of 1. The maximum number of episode time steps is configured at 60 for the first example and 40 for the second example. In addition to these specifications, the Tensorforce agent's default settings are utilized.
Based on our observations for the default choices of hyperparameters in the Tensorforce agent, we refined them step by step to achieve a good balance between exploration and exploitation, thus, leading to our novel agent configuration.

\subsubsection{Software setup}
\label{sec:soft}
\begin{figure}[H]
	\centering
	\includegraphics[scale=.7]{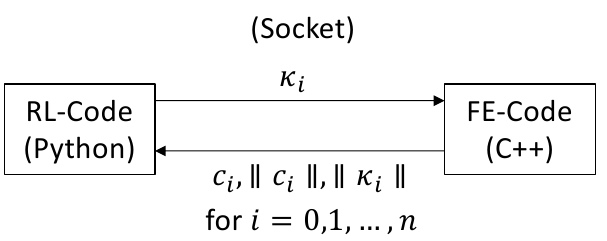}
	\caption{Software setup}
	\label{Fig:Softwaresetup}
\end{figure}

As highlighted in fig.~\ref{Fig:RLAlgOverview}, for this work we implemented a model-based RL approach. The discretization and solution of the transport problems are performed with an in-house C++-based finite element (FE) code. These fully discrete models represent the environment that the agent is interacting with. Therefore, the simulation code and the ML agent have to communicate with each other. This works in the following way: the Python-based RL code sends the suggested actions, i.e. the proposed diffusivities, via sockets to the FE code and the latter sends back the observed state and the $L^2$ norms of the diffusivity and state in each iteration/time instance $i=0,1,\ldots,n$, as illustrated in fig.~\ref{Fig:Softwaresetup}.

The proposed setup preserves the independence of both codes and only requires the development of an interface that, being based on sockets, is very fast. The motivation for this approach is two-fold. First, it takes advantage of the relative strengths of the Python and C++ programming languages. Second, it shows that the idea can be used to control complex models that are implemented in legacy codes, just by programming the corresponding interface.

\section{Results}
\label{sec:results}
In this section, we leverage the model-based reinforcement learning algorithm detailed in Section~\ref{sec:alg} and the corresponding software setup outlined in Section~\ref{sec:soft} to address two examples aligned with the objectives introduced in Section~\ref{sec:cont}. The first objective targets high diffusivities without increasing infections/temperatures, utilizing the reward function~(\ref{eq:R1}). The second objective, guided by the reward function~(\ref{eq:R2}), aims at obtaining low infections/temperatures without decreasing diffusivities. For both case studies the environment is modeled by problem~(\ref{eq:model}) with $\rho=1$.

Our exploration begins with a straightforward heat conduction problem defined on a square, featuring a specified initial condition. Once we establish the efficacy of our approach in this context, we transition to a more complex disease transport problem, defined on the geographical layout of the map of Spain, also with a predefined initial condition. It is important to note that while the results presented in this section are highly promising, they are merely academic, since they have not been calibrated or validated against real-world data.

For both cases, the diffusivity values are initially randomized, and we train for 1000 episodes.

\subsection{Heat conduction within a block with initial condition}
\label{subsec:block}
In this example, we opt for a simple heat conduction problem defined on a square domain with dimensions $1\times 1$. The initial condition is characterized by a uniform temperature of zero across the entire domain, except within a centered circle with a radius of $0.3$, where the temperature is set to 1. This arrangement establishes a distinct heat source at the center of the domain.

Additionally, Neumann boundary conditions are employed, meaning that the heat flux across the boundary, namely $-\kappa\nabla c \cdot n$ is set to zero. Addressing the constrained RL problem posed by the initial boundary value problem in eq. (\ref{eq:model}), we pursue the two specified objectives. The reaction constant $\beta$ is set to $2.5$.

To avoid bias from the distribution, we scale the action space to operate within the range of $[-1,1]$. However, for result evaluation, we rescale the action space back to its original range.
\begin{figure}[htbp]
    \hspace{-.2em}\includegraphics[scale=.18]{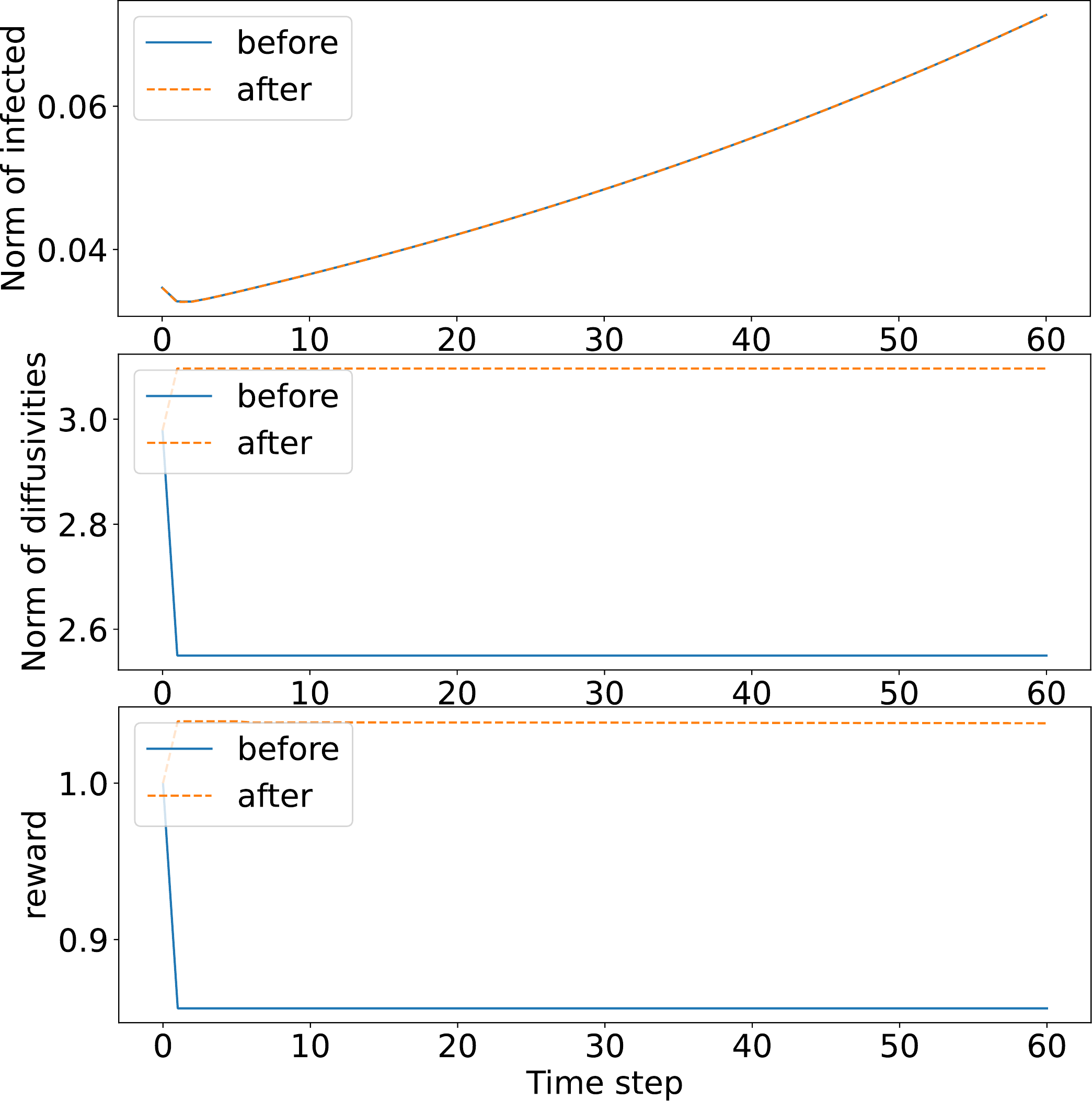}\includegraphics[scale=.18]{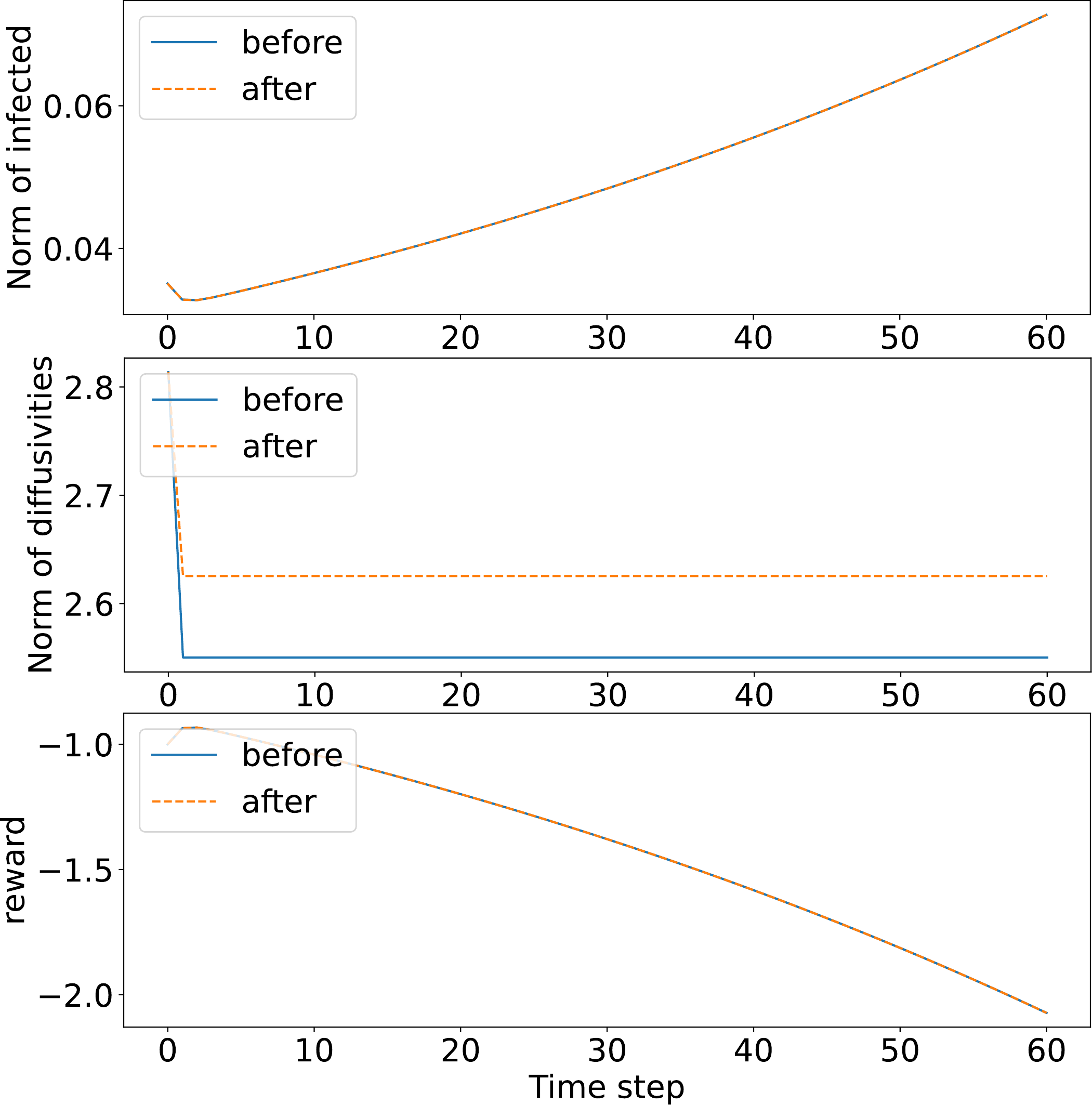}
	\caption{$L^2$ norm of infections and diffusivities and reward function value as in eq. (\ref{eq:R1}) on the left side and as in eq. (\ref{eq:R2}) on the right side from before in orange and after training in blue.}
 \label{Fig:lineplotsIC}
\end{figure}
\begin{figure}[htbp]
    {\centering \includegraphics[scale=.3]{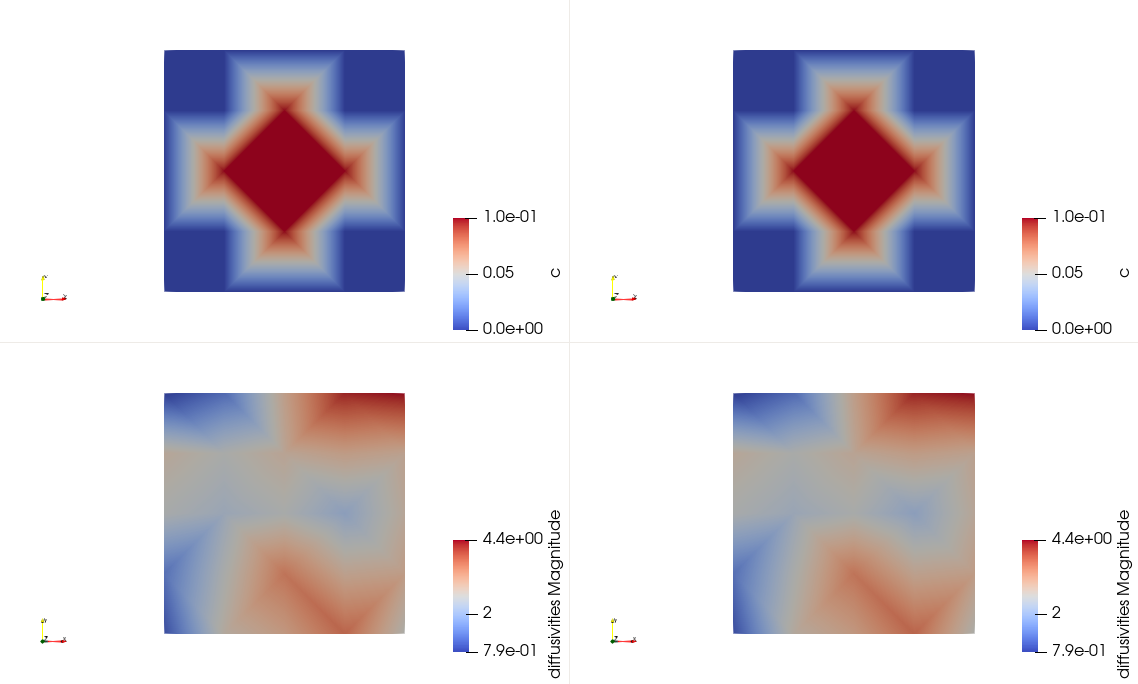}}\\
    time 0
\hrule
    {\centering\includegraphics[scale=.3]{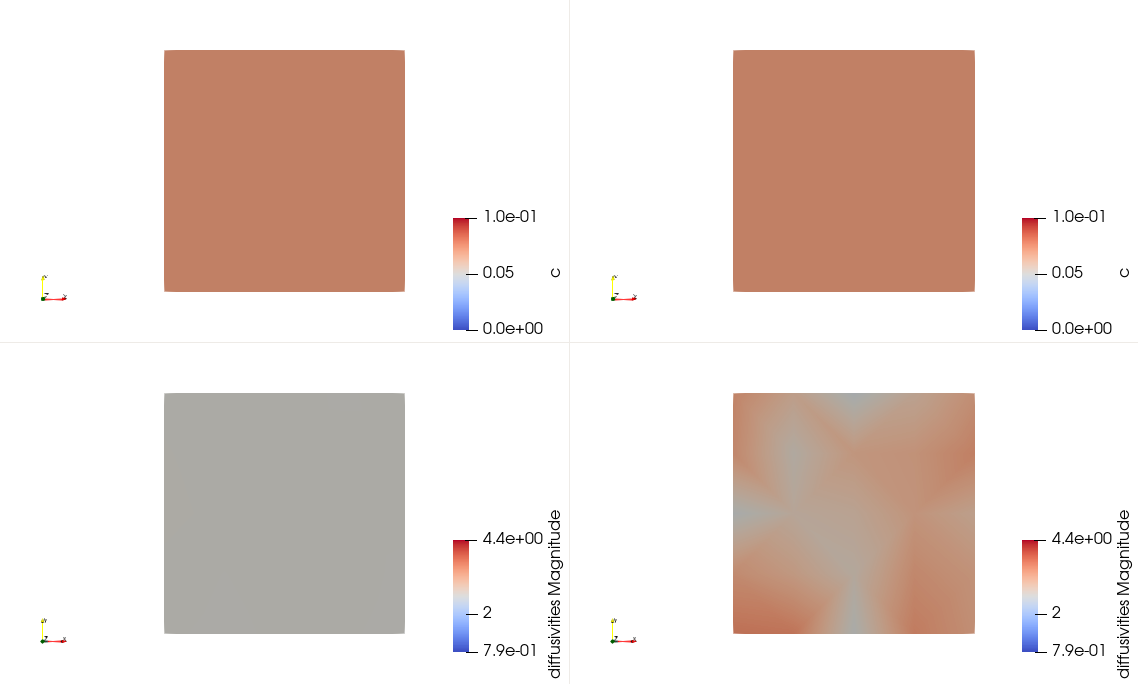}}\\
     time 60
	\caption{RL simulation results for state $c$ and diffusivities at start time on the top and at final time on the bottom resulting from choosing reward function eq. (\ref{eq:R1}).}
  \label{Fig:simresIC}
\end{figure}
\begin{figure}[htbp]
	{\centering\includegraphics[scale=.3]{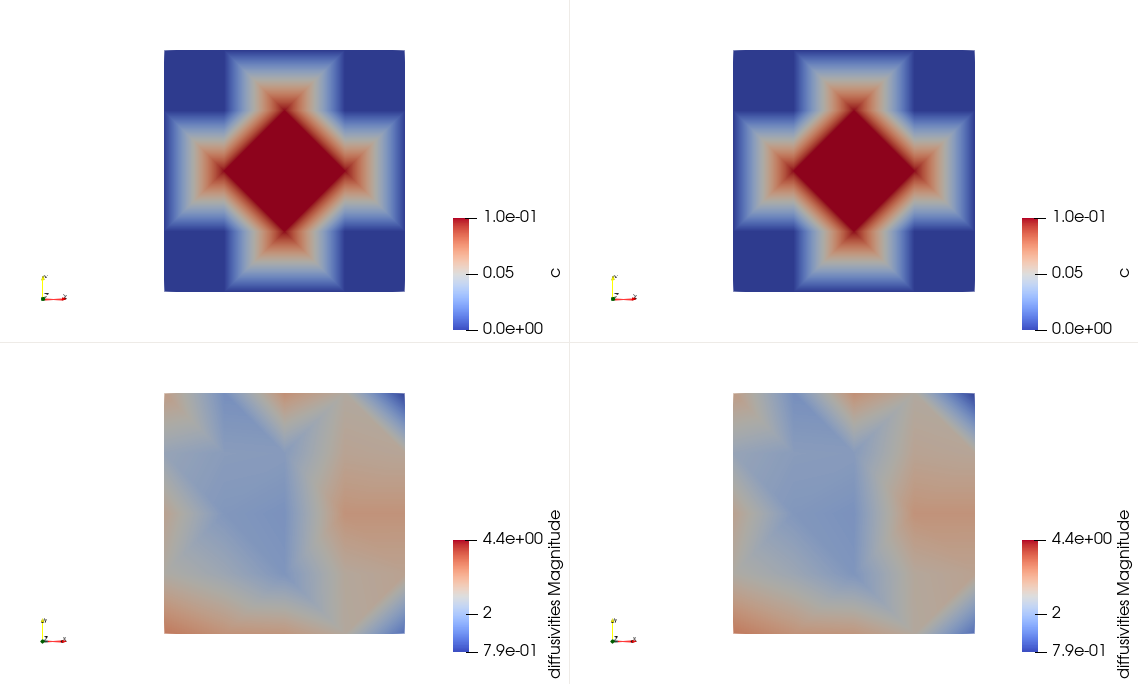}}\\
 time 0
\hrule
    {\centering\includegraphics[scale=.3]{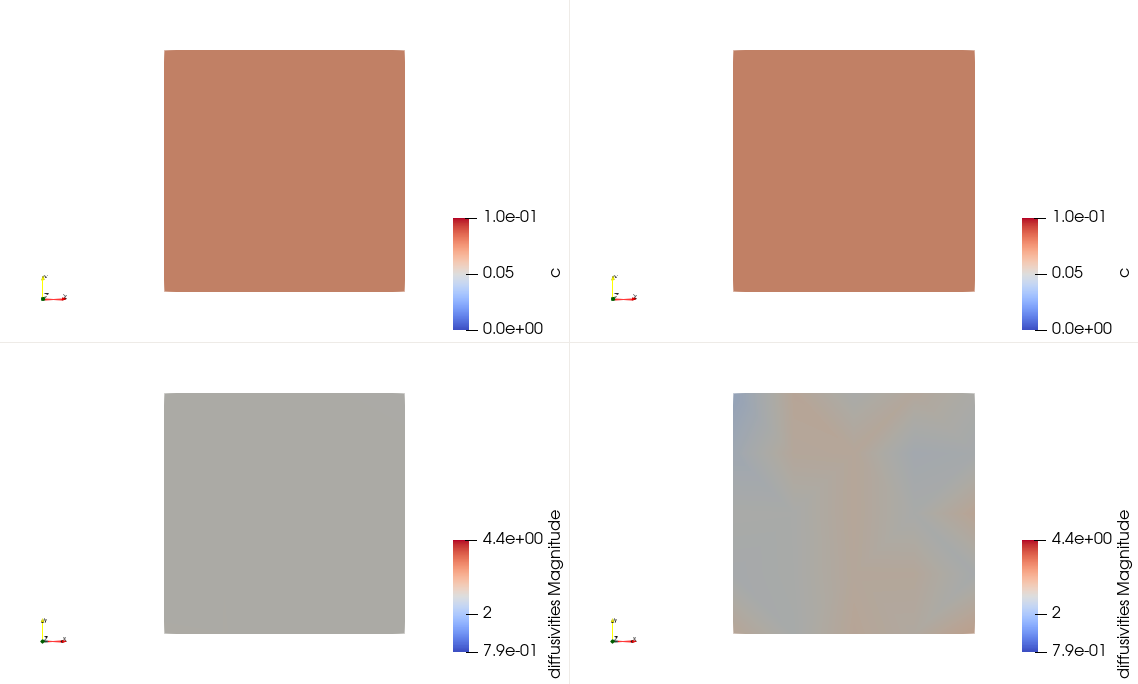}}\\
    time 60
	\caption{RL simulation results for state $c$ and diffusivities at start time on the top and at final time on the bottom resulting from choosing reward function eq. (\ref{eq:R2}).}
  \label{Fig:simresICtheta}
\end{figure}
\begin{figure}[htbp]
    \includegraphics[scale=.22, angle=90]{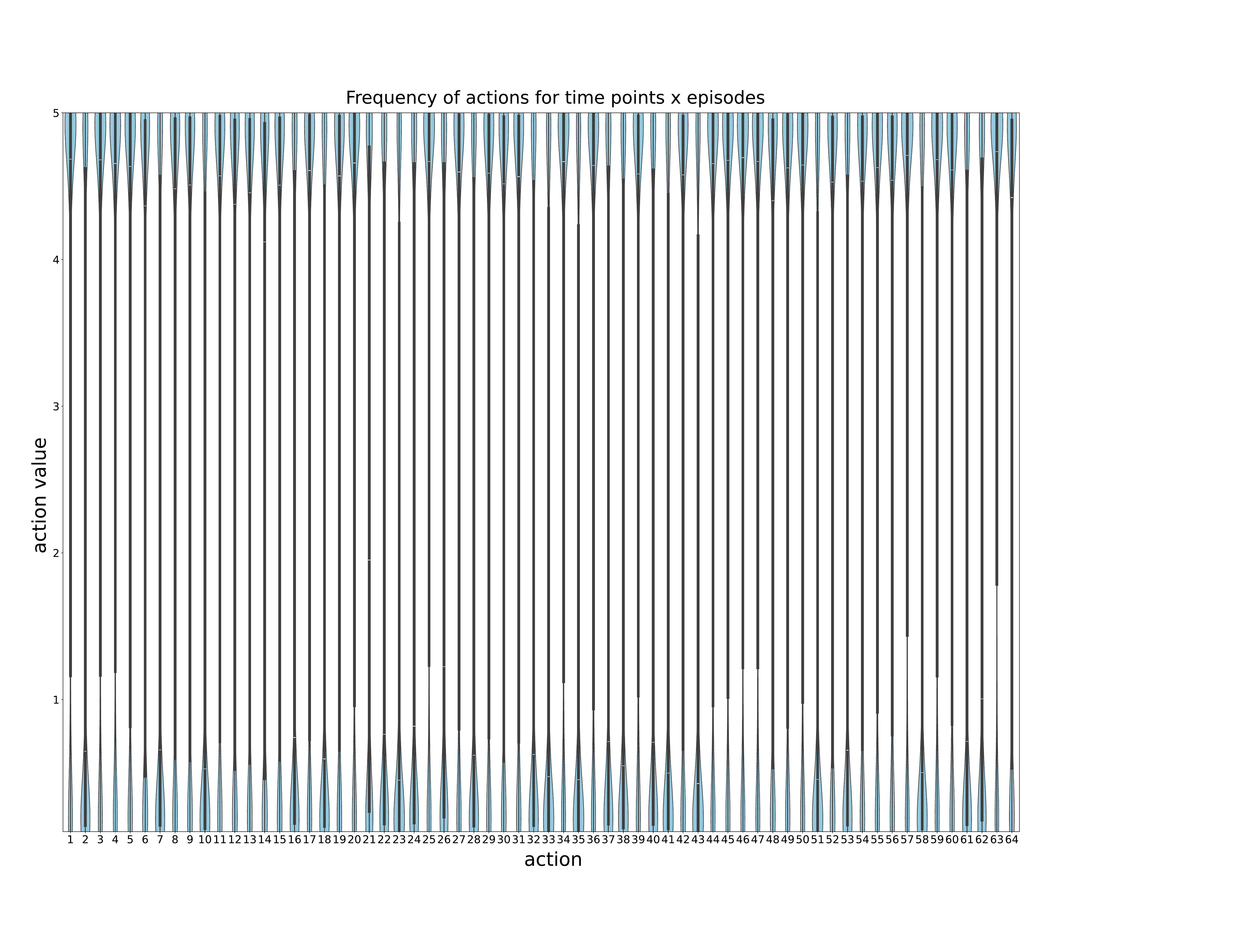}
    \caption{Violin plot showing the distributions of actions for time points times episodes choosing reward function eq. (\ref{eq:R1})}
    \label{fig:violinplotsIC1}
\end{figure}
\begin{figure}[htbp]
    \includegraphics[scale=.22, angle=90]{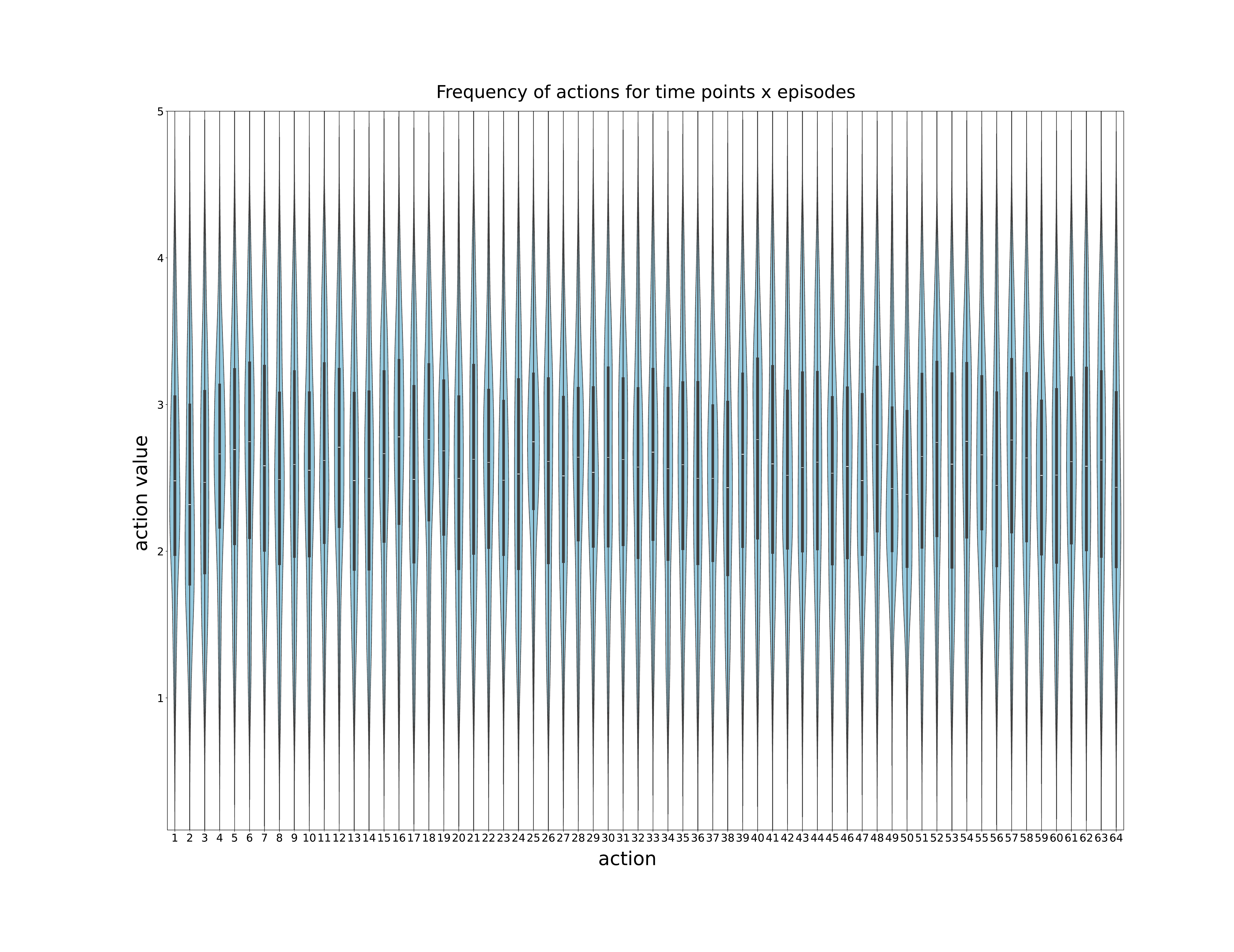}
    \caption{Violin plots showing the distributions of actions for time points times episodes choosing reward function eq. (\ref{eq:R2})}
    \label{fig:violinplotsIC2}
\end{figure}

Fig. ~\ref{Fig:lineplotsIC} shows the $L^2$ norms of the states and diffusivities and the values of the reward function before and after training for the two reward functions, i.e. eq.~(\ref{eq:R1}) on the left and eq.~(\ref{eq:R2}) on the right. When optimizing for mobility, i.e. eq.~(\ref{eq:R1}), diffusivities grow while keeping infections at the same level as before training. Hence the reward function values are significantly higher after training compared to before training. When minimizing the states via eq.~(\ref{eq:R2}), higher diffusivities are possible while keeping the infections as before training.

These observations are confirmed by the simulation results at the start and final time for the two objectives highlighted by fig.~\ref{Fig:simresIC} and \ref{Fig:simresICtheta}. In fig.~\ref{Fig:simresIC} (time 60), before training, the 64 diffusivities are uniformly low across the entire domain, while post-training, they significantly increase across a larger portion of the domain. In fig.~\ref{Fig:simresICtheta} the diffusivities after training are chosen higher as well but significantly lower and in smaller areas of the domain than for the other objective.

The accompanying violin plots (fig.~\ref{fig:violinplotsIC1} and fig.~\ref{fig:violinplotsIC2}) further illuminate these distinctions. For the first objective, two distinct violins emerge, indicating a preference for either low or high diffusivities across all actions. In contrast, the second objective yields a single violin centered within the action space, suggesting a more concentrated and balanced diffusivity distribution.

\subsection{Infection spread in Spain}
\label{subsec:spain}
Transitioning from a simple heat conduction problem, our focus extends to the study of infection spread in Spain. Political decisions during events like the COVID-19 pandemic often operate on a regional level in Spain, with each of the 17 autonomous communities implementing differing measures. To facilitate these decisions, we delve into optimal control policies aiming to strike a balance between economic and healthcare objectives. These policies involve controlling regional diffusivities to minimize infections while maximizing mobility.

To simplify our analysis, we concentrate on the 15 autonomous communities located on the Spanish peninsula, excluding the islands. This selection aligns with the mesh of Spain's map, which consists of these 15 inland regions or autonomous communities as detailed in \cite{mapwiki}: Galicia, Asturias, Cantabria, the Basque Country, Navarra, La Rioja, Aragon, Catalonia, Castilla y Leon, Extremadura, Castile-La Mancha, Madrid, Valencia, Andalusia and Murcia.

\begin{figure}[htbp]
    \centering
    \includegraphics[scale=.125]{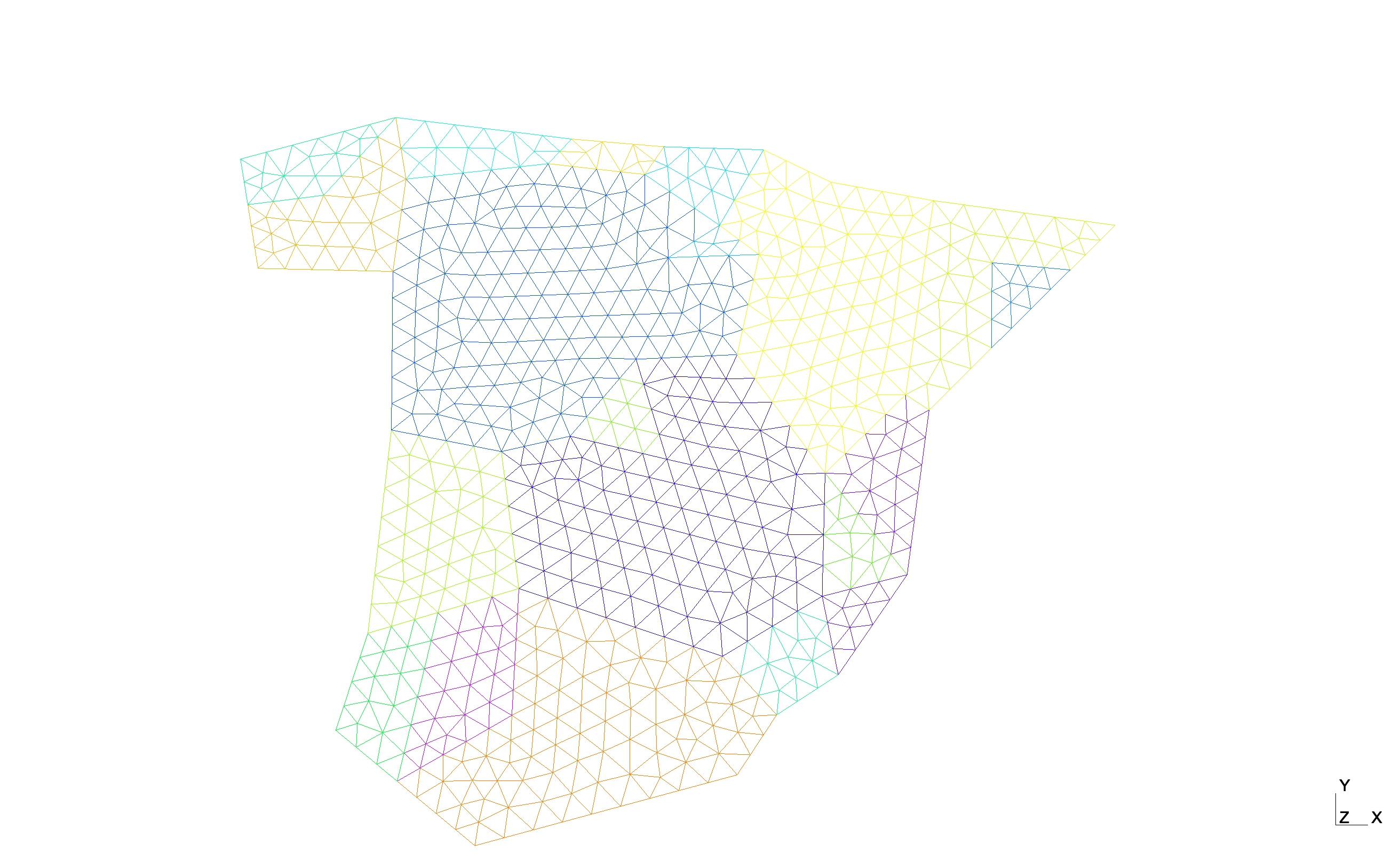}
    \caption{Mesh used for generating simulation results with the provinces highlighted in different colors.}
    \label{fig:mesh}
\end{figure}

This approach diverges from our previous scenario, where the reinforcement learning (RL) algorithm controlled diffusivity values at all mesh points. Instead, we enhance problem-solving tractability for the large-scale problem by having the RL algorithm send only 15 diffusivities to the finite element (FE) code at each time step, corresponding to each autonomous community in Spain. Fig. \ref{fig:mesh} illustrates the mesh used to generate the subsequent results.

As in the previous example, we employ eq. (\ref{eq:model}) to represent the environment within our reinforcement learning algorithm, where we select the contact rate $\beta=50$ based on our experience with this model. In addition, we choose to use an initial condition where there are no contagions anywhere in Spain, except for the central region of Madrid where the infection ratio is $c=1$. In the case of boundary conditions, we assume that there is no flow of infected people entering or leaving the country, so $\nabla c \cdot n =0$. Similar to the previous example in Section~\ref{subsec:block}, we address the initial boundary value problem-based RL problem employing the reward functions defined in eq. (\ref{eq:R1}) and eq. (\ref{eq:R2}). Moreover, to avoid further issues related to the curse of dimensionality, we keep choosing actions from the space $[0.1,5]$ as in the previous example, but here, when handing over the actions to the FE code, we scale them up to the corresponding diffusivity parameters for the map of Spain that are of orders of magnitude 3 to 4, i.e. $\kappa\in[1000, 50000]$, by multiplying them by a constant.
Additionally, we scale the action space to $[-1,1]$ for the algorithm and subsequently rescale it back to the original RL action space of actions in [0.1, 5] for result evaluation.

 With these data, the FE solver can update the value of the infection ratio at every point and send it back to the agent along with its $L^2$ norm and the $L^2$ norm of the diffusivities (see again fig.~\ref{Fig:Softwaresetup}).
\begin{figure}[htbp]
	\centering
    \includegraphics[scale=.24]{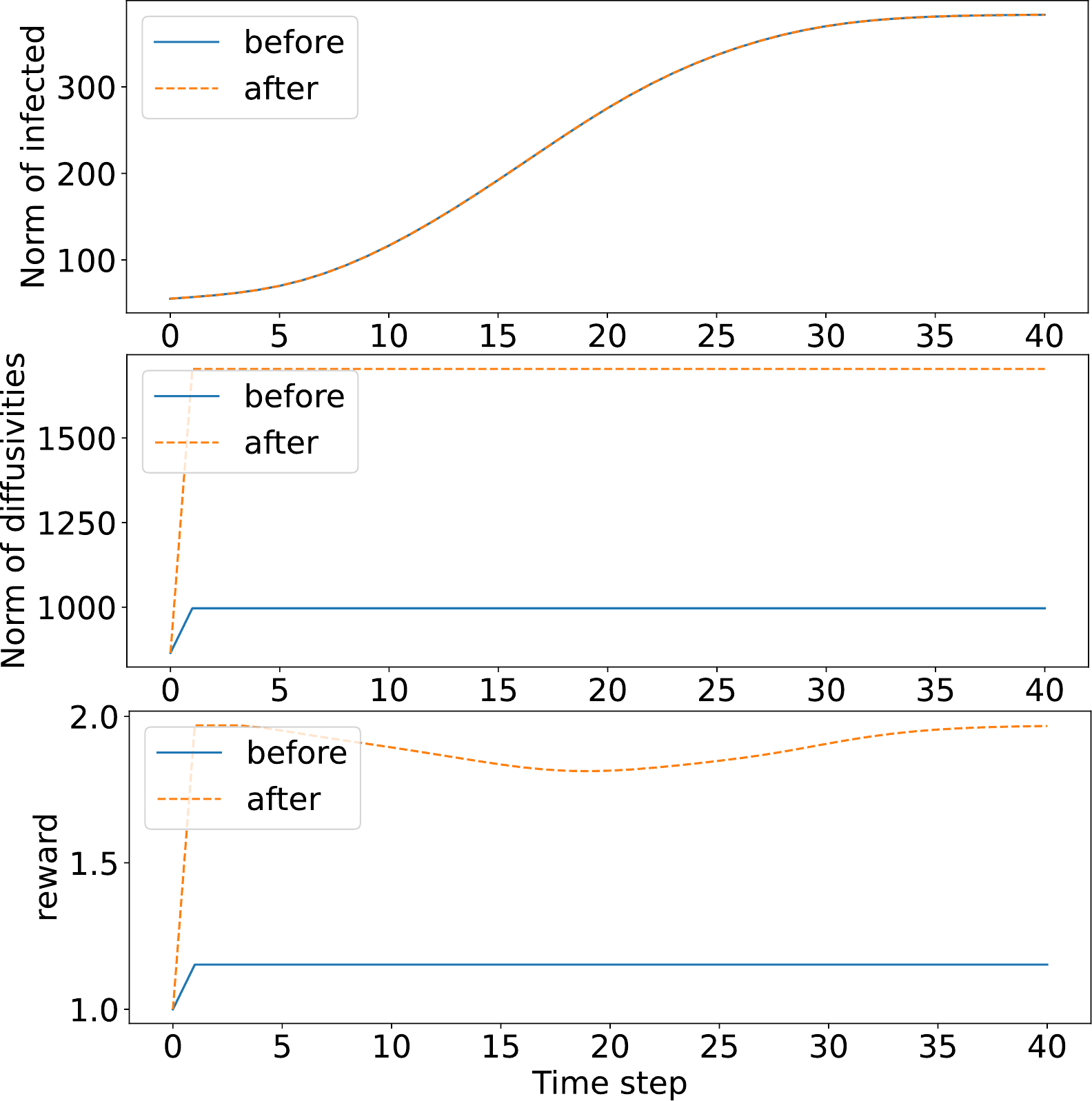}\hspace{.1em}\includegraphics[scale=.24]{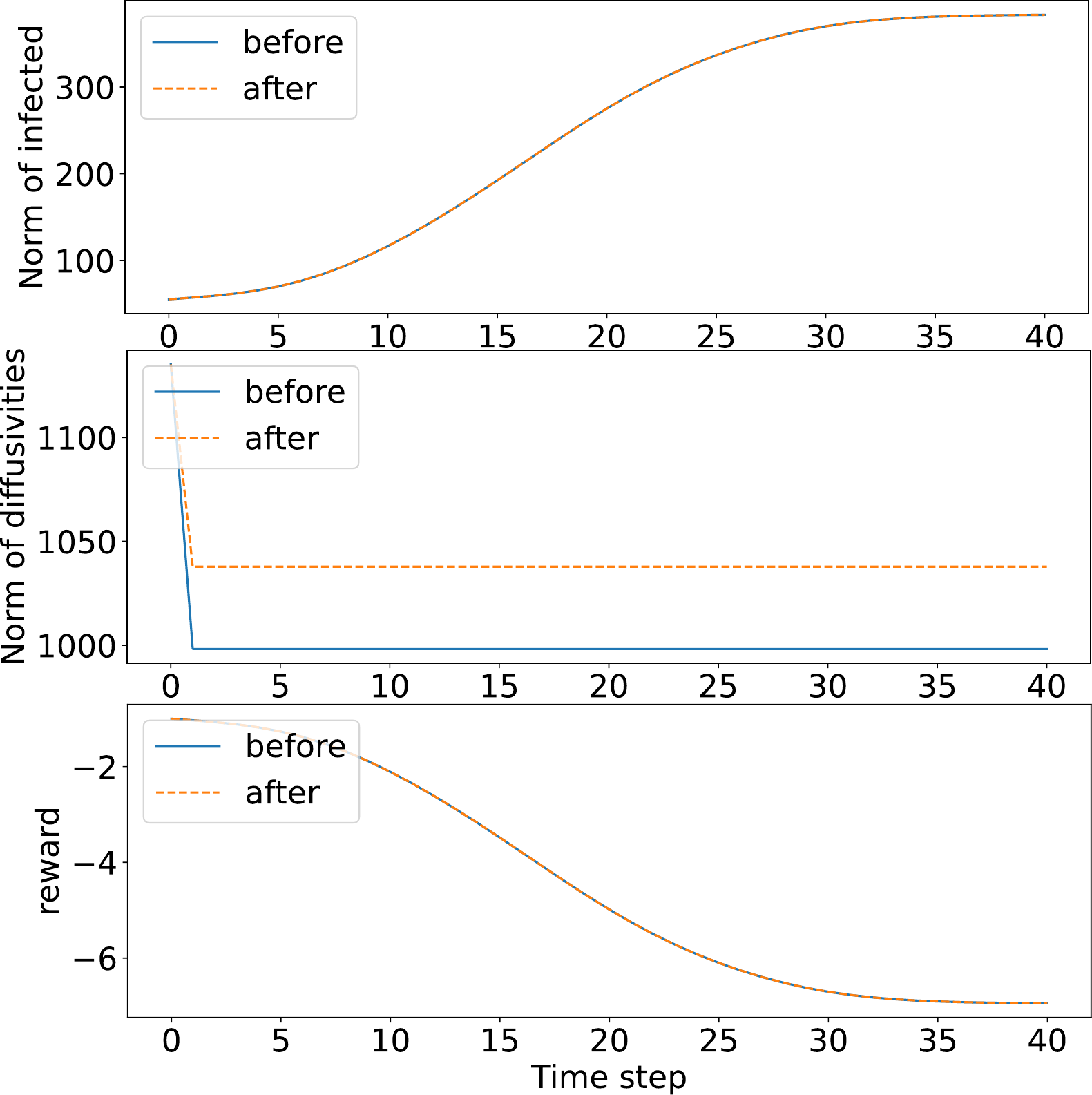}
	\caption{$L^2$ norm of infections, diffusivities, and reward function value as in \ref{eq:R1} on the left side and as in \ref{eq:R2} on the right side from before and after training.}
 \label{Fig:lineplotsSpain}
\end{figure}
\begin{figure}[htbp]
	\centering
    time 0\includegraphics[scale=.13]{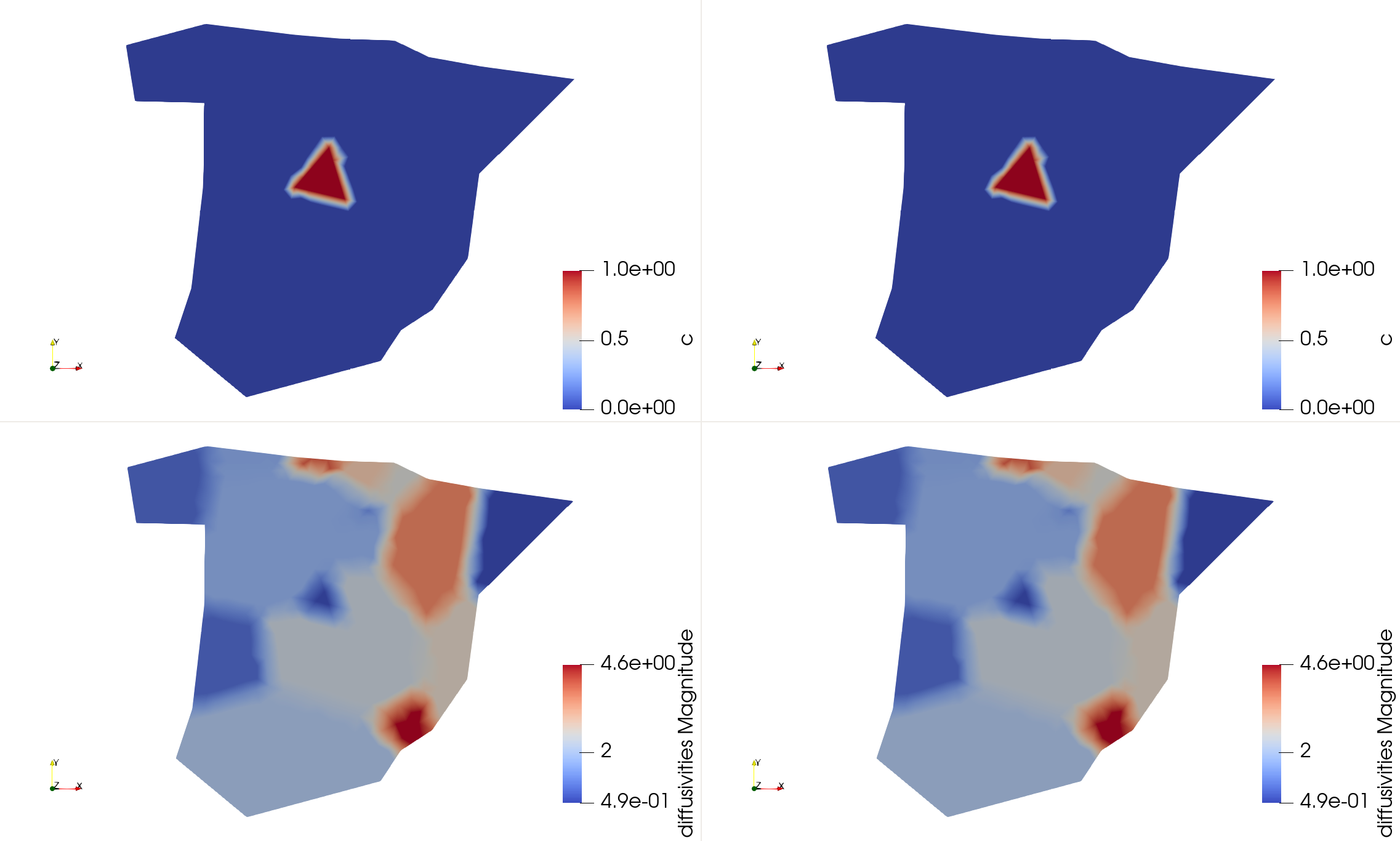}
    \hrule
    time 20\includegraphics[scale=.13]{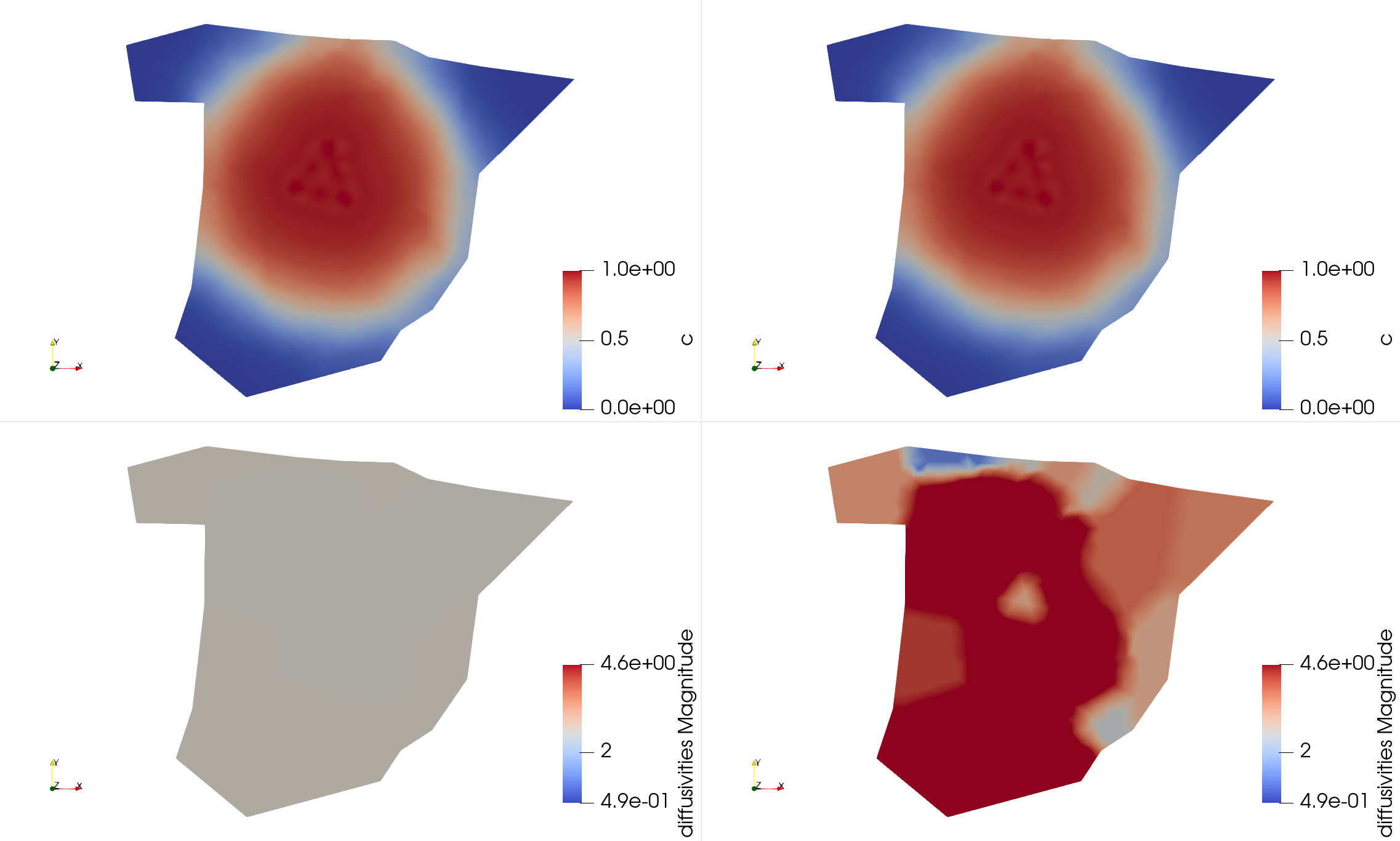}
    \hrule
    time 40\includegraphics[scale=.13]{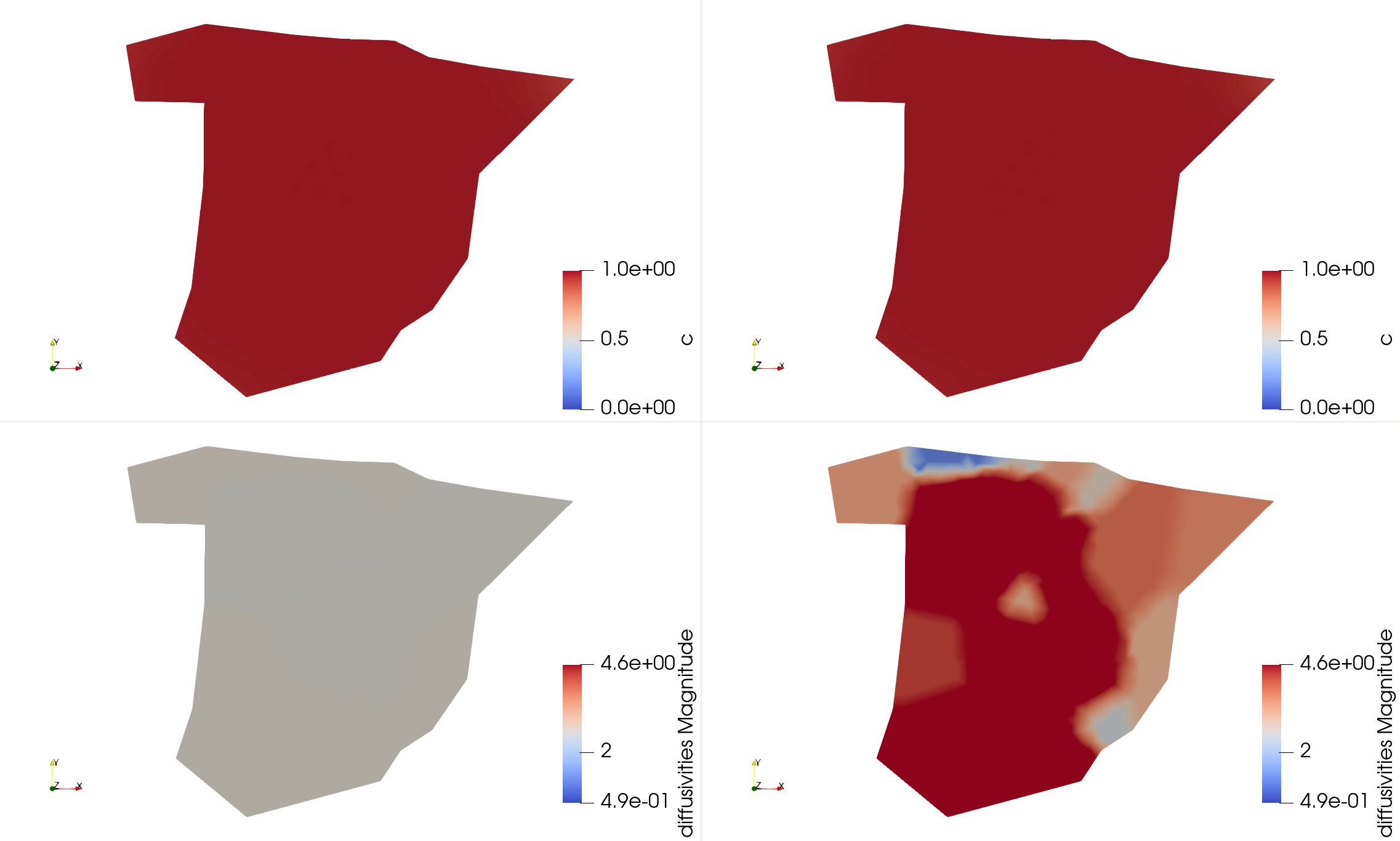}
	\caption{RL simulation results for state $c$ and diffusivities at time 0, 20, and 40 resulting from choosing reward function eq. (\ref{eq:R1})}
 \label{Fig:simresSpain}
\end{figure}
\begin{figure}[htbp]
	\centering
    time 0\includegraphics[scale=.13]{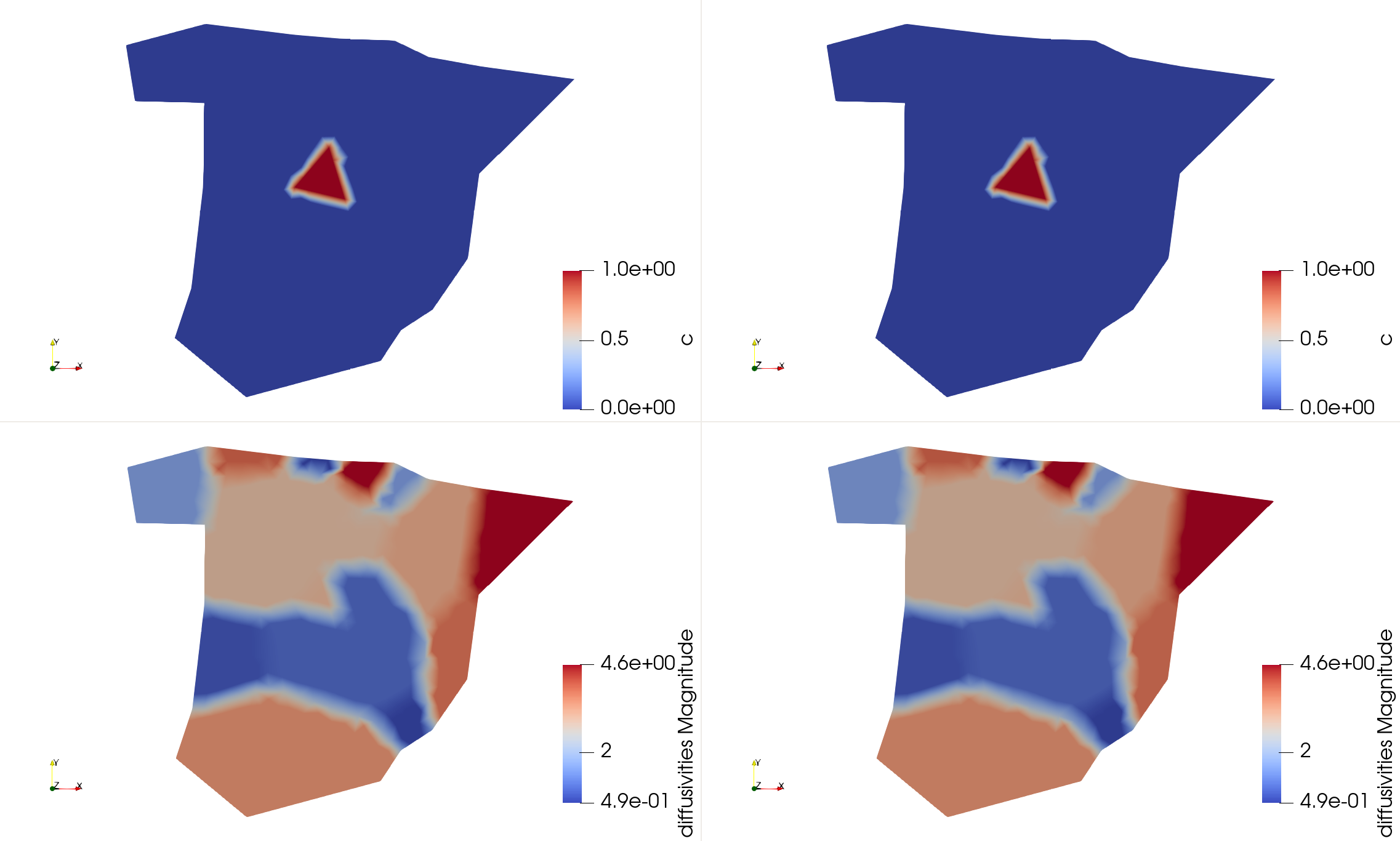}
    \hrule
    time 20\includegraphics[scale=.13]{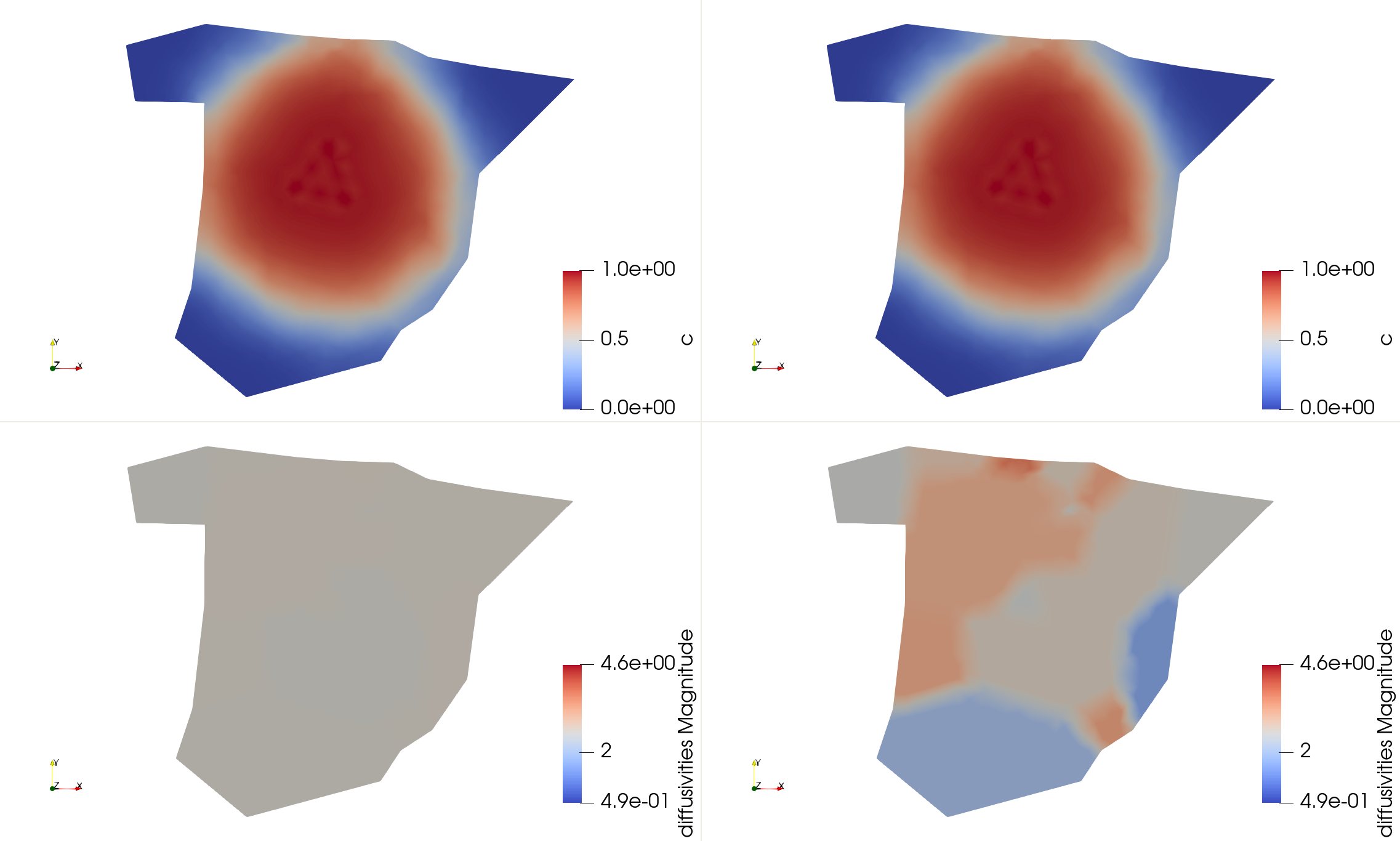}
    \hrule
    time 40\includegraphics[scale=.13]{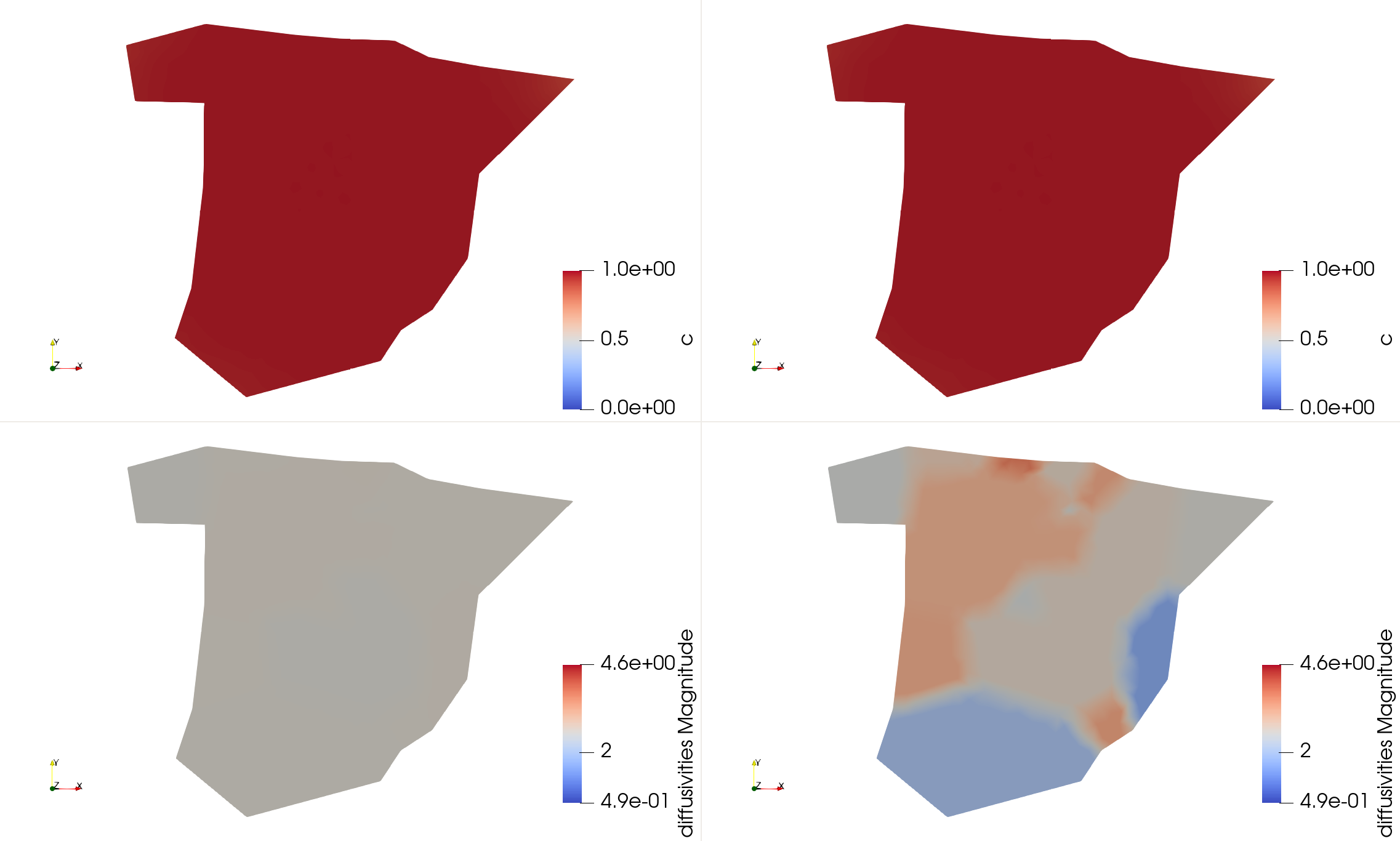}
	\caption{RL simulation results for state $c$ and diffusivities at time 0, 20, and 40 resulting from choosing reward function eq. (\ref{eq:R2})}
 \label{Fig:simresSpaintheta}
\end{figure}
\begin{figure}[htbp]
    \centering
    \includegraphics[scale=.22]{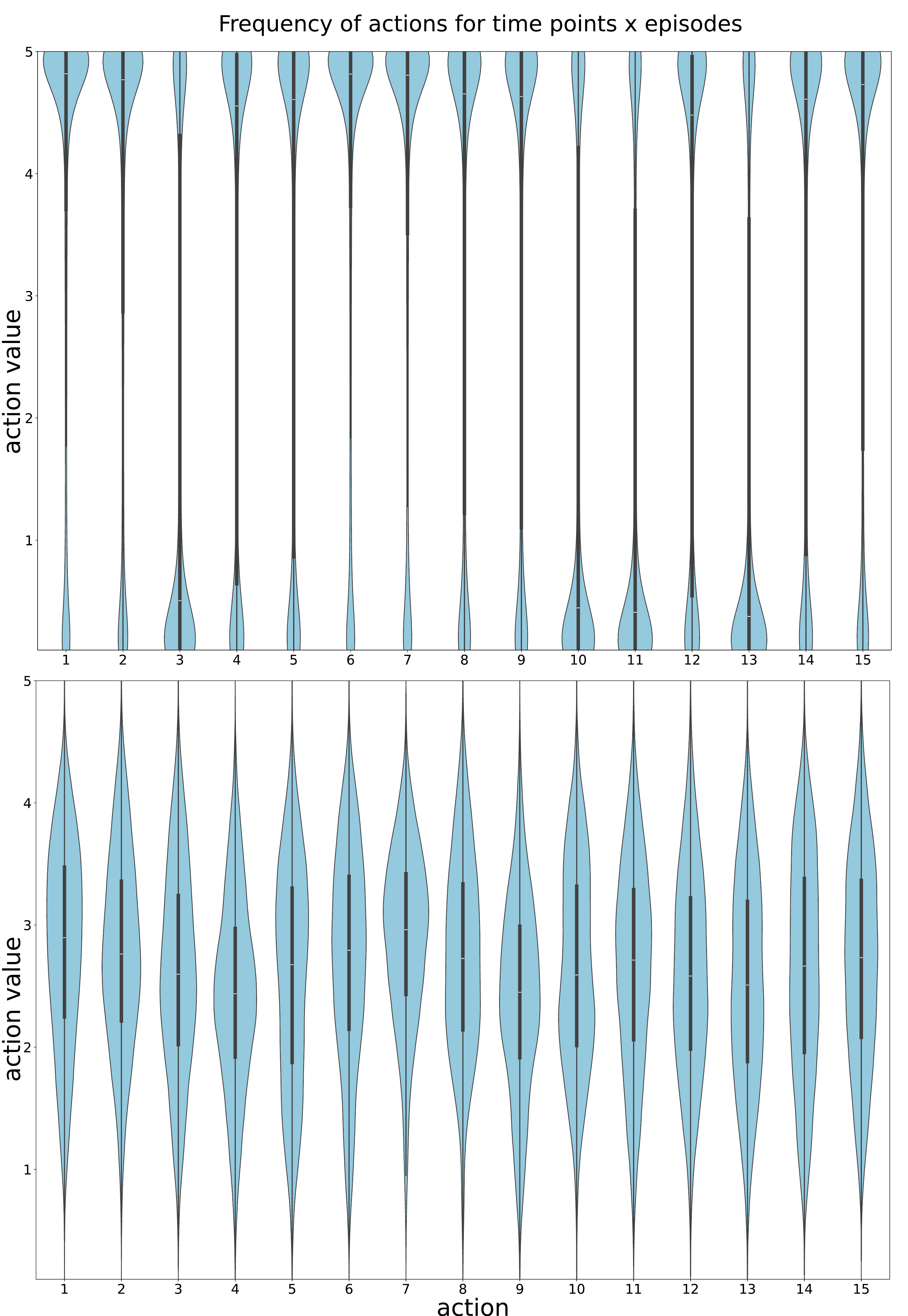}
    \caption{Violin plots showing the distributions of actions for time points times episodes choosing reward function eq. (\ref{eq:R1}) on the top and choosing the reward function eq. (\ref{eq:R2}) on the bottom.}
    \label{fig:violinplotsSpain}
\end{figure}

The plots depicted in fig.~\ref{Fig:lineplotsSpain} show the $L^2$ norms of the states and diffusivities, along with the values of the reward functions for eq. (\ref{eq:R1}) and eq. (\ref{eq:R2}). To maximize mobility and minimize infections, significantly higher diffusivities are selected for eq. (\ref{eq:R1}), and marginally higher diffusivities are chosen for eq. (\ref{eq:R2}), while maintaining infections at a consistent level as observed before training. Consequently, the values of the reward function for eq. (\ref{eq:R1}) exhibit a substantial increase after training compared to their pre-training counterparts.

Fig.~\ref{Fig:simresSpain} and \ref{Fig:simresSpaintheta} validate these observations. Specifically, in fig.~\ref{Fig:simresSpain}, the diffusivities, initially randomly distributed among the 15 regions, are uniformly low throughout the entire domain before training. However, the diffusivities are notably increased post-training across most of the domain, especially in regions with low infection rates. On the contrary, the agent sets higher diffusivities in most regions with high infection rates. This behavior can be explained with the following rationale: In highly infected regions, diffusivities are kept high, since reducing them will not improve the already bad situation, and, at least, the agent attempts to maximize mobility (cf. fig.~\ref{Fig:simresSpain}, time 20). Regions that do not border highly infected areas show low to moderate diffusivities; for example, Asturias, Murcia and the Basque Country exhibit low diffusivities. In fig.~\ref{Fig:simresSpaintheta}, when optimizing for low infections, after training, regions that are already infected have high or marginal diffusivities and some regions that are not infected yet, such as Andalusia and Valencia show low diffusivities (cf. fig. \ref{Fig:simresSpaintheta}, time 20).

The distributions for the 15 actions, i.e. the diffusivities for each of the 15 regions, highlighted by the violin plots fig.~\ref{fig:violinplotsSpain} show, similar to the previous example from Section~\ref{subsec:block} two violins for the first objective, meaning that either high or low diffusivities are favored for all 15 regions. Conversely, for the second objective, the algorithm prefers medium diffusivities, illustrated by a much broader violin across the domain in fig.~\ref{fig:violinplotsSpain}.

\section{Discussion}
\label{sec:concl}
In this article, we studied the application of reinforcement learning control to transient diffusion-reaction models. The analysis presented addresses them from various perspectives such as mathematical modeling, discretization, algorithmic considerations, and software implementation.

After introducing the mathematical modeling, discretization techniques, and control formulations, we delved into the specifics of the RL algorithm settings and the intricacies of the software setup designed to facilitate communication between Python and C++-based codes. The primary outcome of this research lies in the exploration of a novel RL agent, achieved through a modification of the standard Tensorforce agent. The latter, coupled with new reward functions, showcased the successful application of RL-based control to address thermal and disease transmission dynamics problems.

The results of this study suggest that RL has significant potential for application in various fields, particularly in aiding decision-making processes. However, as demonstrated in this work, simplifications in the modeling phase may be required to enhance the tractability of the solution space. Additionally, the findings underscore the ongoing need for advancements in algorithm reliability and explainability to further solidify the applicability of RL across various domains.
\backmatter
\bmhead{Acknowledgements}
The authors acknowledge the financial support from the Madrid regional government through a grant from IMDEA Materials that addressed research activities on SARS-CoV-2 and COVID-19 and was financed with REACT-EU resources from the European regional development fund.

\bibliography{cas-refs}

\end{document}